\documentclass[12pt]{amsart}

\usepackage{amsmath,amssymb,amscd}
\input xy
\xyoption{all}

\theoremstyle{plain}
\newtheorem{theorem}{Theorem}[section]

\newtheorem{lem}[theorem]{Lemma}

\newtheorem{prop}[theorem]{Proposition}

\newtheorem{cor}[theorem]{Corollary}
\theoremstyle{definition}

\newtheorem{rem}[theorem]{Remark}

\newcommand{\ra}{\rightarrow}

\DeclareMathOperator{\Fix}{{Fix}}

\DeclareMathOperator{\Ind}{{Ind}}

\DeclareMathOperator{\Div}{{Div}}
\DeclareMathOperator{\End}{{End}}
\DeclareMathOperator{\rk}{{rk}}

\newcommand{\R}{\mathbb R}

\def\Z{{\mathbb Z}}
\def\Q{{\mathbb Q}}

\def\P{{\mathbb P}}
\def\C{{\mathbb C}}

\def\A{{\mathcal A}}  

\def\A{{\mathcal A}} 
\def\bs{{\backslash}}

\numberwithin{equation}{section}

\begin{document}
\title{A Galois-theoretic approach to Kanev's correspondence}

\author{Herbert Lange}

\address{Mathematisches Institut, Bismarckstr. 1 1/2, 91054 Erlangen, Germany}

\email{lange@mi.uni-erlangen.de}

\author{Anita M. Rojas}

\address{Departamento de Matematicas, Facultad de Ciencias
Universidad de Chile}

\email{anirojas@uchile.cl} 

\thanks{This work was supported by FONDECYT No. 11060468}

\begin{abstract}
Let $G$ be a finite group, $\Lambda$ an absolutely irreducible $\Z[G]$-module and $w$ a weight of $\Lambda$. 
To any Galois covering with group $G$ we associate two correspondences, the Schur and the Kanev correspondence. 
We work out their relation and compute their invariants.
Using this, we give some new examples of Prym-Tyurin varieties.

\end{abstract}


\subjclass[2000]{Primary: 14K05; Secondary: 14H40}

\maketitle

\section{Introduction}

Recall that a {\it Prym-Tyurin variety of exponent} $e$ in the Jacobian $JC$ of a smooth projective curve $C$ is by definition 
an abelian subvariety $P$ of $JC$ such that the canonical principal polarization restricts to the $e$-fold of a 
principal polarization on $P$.
Starting from a complex semisimple Lie algebra $\bf g$, 
a rational map $f: \C \ra \bf g$ such that all but a finite number of elements $f(\xi)$ 
are regular and semisimple, and a faithful irreducible representation 
$\rho: {\bf g} \ra {\bf gl}(V)$, Kanev associated in \cite{kan}
a spectral cover $C \ra \P^1$. Moreover, denoting by $W({\bf g})$ the Weyl group of $\bf g$, given a $\Z[W({\bf g})]$-module
$\Lambda$ induced by the representation $\rho$, he constructed a correspondence on the curve $C$ which under suitable 
assumptions defines a Prym-Tyurin variety in the Jacobian $JC$. \\

In this paper we consider the Kanev correspondence from a slightly different point of view. Let $G$ be an arbitrary finite group. 
We start with a Galois covering $\pi: X \ra \P^1$ with Galois group $G$ and an absolutely irreducible $\Z[G]$-module $\Lambda$.
A weight of $\Lambda$ is by definition an element $w \in W = \Lambda \otimes_{\Z} \Q$ satisfying $gw-w \in \Lambda$ for all 
$g \in G$. In the case of a Weyl group this is equivalent to the usual definition.

To every such triple $(G,\Lambda,w)$ one can associate two correspondences in a natural way. First it defines a 
canonical negative definite $G$-invariant symmetric bilinear form $(\;,\;)$. Schur's orthogonality relations 
induce a correpondence $S_w$ on the curve $X$. If $H$ denotes the stabilizer subgroup of $w$ in $G$, 
we call the push-down $\overline{S}_w$ of $S_w$ to the quotient curve $C = X/H$ the {\it Schur correspondence} associated to the triple.     
On the other hand, Kanev's construction gives a correspondence $\overline{K}_w$ on $C$, which we call the {\it Kanev correspondence}
associated to the triple. The main result of the paper is the following relation between both correspondences 
(see Theorem \ref{kanev-schur})
$$ 
\overline{S}_w + |H|^2 (\Delta-\overline{K}_w)=[(w,w)+1]|H|^2 T.
$$
Here $\Delta$ denotes the identity correspondence and $T$ the trace correspondence of the covering $C \ra \P^1$. Noting that  
every correpondence on $C$ defines an abelian subvariety of the Jacobian $JC$, we obtain as a consequence that the the abelian subvariety 
associated to the Schur and Kanev correspondences coincide. Moreover we apply this relation in order to compute the invariants of 
the Kanev correspondence.\\

In the last part of the paper we use our results and Kanev's criterion (see \cite{K1}) 
in order to construct new families of Prym-Tyurin varieties. 

In the first set of examples $G$ is the Weyl group of type $E_6$ 
which was already considered in \cite{kan}, where every covering $C \ra \P^1$ with 
such a monodromy group with ''simplest ramification'' provides 
a Prym-Tyurin variety. We obtain other examples, since the formula of Proposition \ref{prop4.2} allows us to admit more general 
coverings (see Remark \ref{rem5.4}).

Another family of examples uses the alternating group $A_n,\; n \geq 4$ which is not a Weyl group. Although the representation
$\Lambda$ is the restriction of an absolutely irreducible representation of the symmetric group $S_n$, the 
Prym-Tyurin varieties cannot be constructed using this group, since first of all the coverings are different and moreover the 
ramification is not of simplest type. Instead, the ramification is given by 3-cycles. \\

The contents of the paper is as follows: In Section \ref{S:schur} we introduce the Schur correspondence and prove a formula which will be used to 
prove the relation with the Kanev correspondence. Section 3 contains the definition of the Kanev correspondence and the proof of
the above mentioned relation. In Section 4 we compute the invariants of the correspondences. Finally Section 5 
contains the examples.

\section{The Schur correspondence}\label{S:schur}

Let $G$ be a finite group and $\Lambda$ an absolutely irreducible $\Z$-representation of rank $n$ of $G$.
So $\Lambda$ is a lattice on which there is a left action of $G$: $(g,\lambda) \mapsto
g \cdot \lambda = g \lambda$. The action extends to an action of $G$ on the $\Q$-vector space
$$
W=\Lambda \otimes_{\Z} \Q .
$$
The elements $w\in W$ with the property $gw-w \in \Lambda$ for all $g \in G$ are called
the {\it weights} of the representation $\Lambda$. The weights of $\Lambda$ form a lattice $P$ containing
$\Lambda$ with $P/\Lambda \simeq H^1(G,\Lambda)$. In the special case that $G$ is the Weyl group of a root system,
this notion coincides with the usual notion of weights (see \cite{kan}, Prop. 3.2).
Since $W$ is an absolutely irreducible $\Q$-representation, there is, up to a constant, a unique
$G$-invariant bilinear form on $W$.

Now let $w\in W$ be a weight, which is fixed in the sequel. Then there is a unique negative definite $G$-invariant
symmetric bilinear form $(\;,\;)$ on $W$ such that

\begin{enumerate}
\item $(w,\lambda) \in \Z$ for all $\lambda \in \Lambda$, and
\item any $G$-invariant symmetric bilinear form on $W$ with property (1) is an integer multiple of $(\,,)$.
\end{enumerate}

One can associate to any weight $w$ a projector of the rational representation ring $\Q[G]$ by
$$
p_w=\frac{\rk \Lambda}{|G|(w,w)}\sum_{g \in G} (w,gw)g.
$$
The fact that $p_w$ is a projector, is equivalent to Schur's orthogonality relations (see \cite{serre} Ch.2, Thm. 8).\\

Consider a Galois covering
$$\pi:X \to \P^1$$
with Galois group $G$. Thus $X$ is a smooth projective curve with a left action of $G$,
denoted by $(g,x) \mapsto g(x)$, with quotient $\P^1$. The projector $p_w$ of $\Q[G]$ induces
a correspondence with rational coefficients on $X$ over $\P^1$, namely
$$
S_w=\sum_{g \in G} (w,gw)\Gamma_g.
$$
Here $\Gamma_g \in X \times X$ denotes the graph of the automorphism $g$ of $X$.
Considered as a map of $X$ into the group $Div(X)$ of divisors of $X$, $S_w$ is given by
$$
S_w(x)=\sum_{g \in G}(w,gw)g(x).
$$
The $G$-invariance of $(\,,)$ implies that $S$ is symmetric. Moreover
note that $\sum_{g \in G}gw$ is $G$-invariant and thus 0 in $W$, which implies
$$
\deg S_w = \sum_{g \in G}(w,gw)=(w,\sum_{g \in G}gw)=0.
$$
Because of its connection to Schur's relations
we call $S_w$ the {\it Schur correspondence} of $X$ associated to the weight $w$.
It is a correspondence with rational coefficients.\\
Consider the subgroup $H=Stab_G(w)$ of $G$. Denoting $C=X/H$ we
get a diagram


\begin{equation} \label{eq:diagram} \xymatrix{ X \ar[dd]_{\pi}
\ar[dr]^{\varphi} \\ & C \ar[dl]^{\psi} \\ \P^1 }
\end{equation}


\begin{rem}\label{R:sevin}

Recall that the core of a subgroup $G'$ of $G$ is defined as the intersection of all conjugate subgroups of $G'$. 
 With the above notation we have the following facts:\\
(1) $\pi:X \to \P^1$ is the Galois extension of $\psi:C\to \P^1$ if and only if 
Core$\,(H)=1$.\\
(2) Core$\,(H)=\ker(\rho)$, where $\rho:G \to GL(W)$ denotes the homomorphism defined by $W$.\\
This implies that in the case of a faithful representation $W$ the covering $\pi:X\to \P^1$
is the Galois extension of $\psi:C=X/H\to \P^1$. In particular the branch locus of $\pi$ and $\psi$
is the same, and the monodromy group of $\psi$ coincides with the Galois group of $\pi$.

\end{rem}


The correspondence $S_w$ on $X$ descends to a correspondence $\overline{S}_w$ on $C$, namely
$$
\overline{S}_w=(\varphi \times \varphi)_*S_w \subset C\times C.
$$
In order to express this correspondence as a map of $C$ into a symmetric product of $C$,
denote $[G:H]=d$ and let $\{g_1, \dots , g_d\}$ denote a set of representatives for the left cosets of $H$ in $G$:
$$
G= \cup_{i=1}^d g_i H.
$$
\begin{prop}\label{P:schur}
For any $c=\varphi(x) \in C$
$$
\overline{S}_w(c)=|H|^2\sum_{i=1}^d (w,g_iw)\varphi(g_i^{-1}(x)).
$$
\end{prop}
\begin{proof}
Since $H$ is the stabilizer of $w$ we have by definition of $\overline{S}_w$

\begin{eqnarray*}
\overline{S}_w(c) & = & \displaystyle{\sum_{i=1}^d\sum_{h\in H}(w,g_ihw)\sum_{h' \in H}\varphi(g_ihh'(x))}\\
& = & \displaystyle{|H|\sum_{i=1}^d\sum_{k\in H}(w,g_iw)\varphi(g_ik(x))} \\
\end{eqnarray*}

Now, if $\{g_i\}$ is a set of representatives for the left cosets of $H$, then $\{g_i^{-1}\}$
is a set of representatives for the right cosets of $H$. Hence for any pair $(k,i) \in H\times \{1, \dots, d\}$
there is a unique pair $(h,j) \in H \times \{1, \dots, d\}$ such that $g_ik=hg_j^{-1}$.
Moreover, if $g_ik$ runs exactly once through $G$, so do the elements $hg_j^{-1}$. This implies, since
$(\,, )$ is $G$-invariant and since $H$ stabilizes $w$,

\begin{eqnarray*}
\overline{S}_w(c)& = &|H|\sum_{i=1}^d\sum_{h \in H}(w,g_iw)\varphi(hg_i^{-1}(x)) \\
& = &|H|\sum_{i=1}^d\sum_{h\in H}(w,g_iw) \varphi(g_i^{-1}(x))\\
\end{eqnarray*}
which implies the assertion.
\end{proof}

The correspondence $S_{w}$ has degree 0 and thus is not effective. We will use the following lemma to find an
effective correspondence which is equivalent to $S_{w}$ i.e., it induces the same endomorphism $\sigma_{w}$
of the Jacobian $JX$.

\begin{lem}\label{L:Deffective}
{\em (1)} $(w,gw)-(w,w)$ is an integer for all $g\in G$.\\
{\em (2)} $(w,gw)\geq (w,w)$. If $g \in G \setminus H$, then $(w,gw) > (w,w).$\\
\end{lem}

\begin{proof}
(1) is clear, since $(w,gw)-(w,w)=(w,gw-w)\in \Z$ as $w$ is a weight.
For the proof of (2), note
\begin{equation} \label{eq2.1}
0\geq (gw-w,gw-w)=(gw,gw)-2(w,gw)+(w,w)=2[(w,w)-(w,gw)]
\end{equation}
implying $(w,gw)\geq (w,w)$ for all $g \in G$. If $g \in G\setminus H$, we have
$0 > 2[(w,w)-(w,gw)]$ and thus $(w,gw) > (w,w)$.
\end{proof}

Now define a correspondence $D_{w}$ on $X$ by
$$
D_w := S_w - \sum_{g \in G}(w,w)\Gamma_g,
$$
Considered as a map $X \ra Div(X)$ the correspondence $D_w$ is given by
$$
D_w(x):=\sum_{g \in G}[(w,gw)-(w,w)]g(x).
$$

\begin{prop}\label{P:shureffective}
{\em (1)} $D_w$ is an effective symmetric correspondence of $X$ in $X$ of degree $-|G|(w,w)$ with integer coefficients.\\
{\em (2)} The correspondences $S_{w}$ and $D_{w}$ are equivalent.
\end{prop}

\noindent
\begin{proof} According to Lemma \ref{L:Deffective}, we can write
$$
D_{w}=\{(x,y) \in X \times X : y=g(x), (w,gw) > (w,w)\}
$$
and thus $D_{w}$ is effective with integer coefficients.
Next we claim that $D_w$ is symmetric. To see this, let $y=g(x)$. Then
$D_w(y)= \sum_{h \in G}[(w,hw)-(w,w)]h(y)$.
The $G$-invariance of $(\; , \; )$ gives
$$
(w,g^{-1}w)-(w,w)=(w,gw)-(w,w).
$$
Using Lemma \ref{L:Deffective}, this implies $x \in D_{w}(y)$ if and only if $g \notin H$ if and only if $y \in D_{w}(x)$.\\

The degree of $D_{w}$ is given by
$$
\deg D_{w}=\sum_{g \in G}[(w,gw)-(w,w)]=(\sum_{g \in G}gw,w)-|G|(w,w) = -|G|(w,w).
$$
This completes the proof of (1).
For (2) it suffices to note that $\sum_{g \in G} g(x)$ is the pull-back of the point $\pi(x)$ in $\P^1$.
\end{proof}

\begin{rem} \label{rem2.4}
Since the degree of $D$ is large in general, in most cases the correspondence $D_{w}$ will have fixed points.
But there might be examples, where $(w,w)$ is a small positive rational number, in which case the degree is not as large.
\end{rem}

\section{The Kanev correspondence}\label{S:kan}

Let the notation be as in the previous section. We want to study the polarization on the abelian subvariety $B_w$ of $JX$,
associated to a Schur correpondence, which is induced by the principal polarization $\Theta$ of $JX$.
For this, the correspondences $S_w$ and $D_{w}$ seem to be of little use (see Remark \ref{rem2.4}).
Kanev introduced in \cite{kan} an effective integral correspondence on the curve $C$,
which we discuss in this section.\\

As in the last section, let $\{g_1 = 1, g_2, \dots , g_d\}$ denote a set of representatives of the left cosets of $H$ in $G$.

\begin{lem}\label{L:positivness}
The number $(g_iw,g_jw)-(w,w)-1$ is a non-negative integer for all $i\neq j$.
\end{lem}

\begin{proof}
Since the scalar product $(\,,)$ is $G-$invariant, it is sufficient to show that $(w,gw)-(w,w)-1$
is a non-negative integer for all $g \in G\setminus H$. But if $g \in G\setminus H$,
the negative definiteness of $(\;,)$ and (\ref{eq2.1}) imply
$(w,gw) > (w,w)$. On the other hand, $gw-w \in \Lambda$ and then
$(w,gw)-(w,w)-1$
is a non-negative integer for all $g \in G\setminus H$.
\end{proof}

Denote by $U$ the complement of the branch locus $B$ of $\psi$ (or $\pi$) in $\P^1$.
Fix a point $\xi_0 \in U$. Since $H$ is the stabilizer of $w$,
the group $G$ acts on the set $\{w =g_1w,g_2w, \ldots,g_dw\}$ as well as on the fibre $\psi^{-1}(\xi_0)$.
Choosing an element in the fibre
$\psi^{-1}(\xi_0)$ induces a $G$-equivariant bijection
$$
\{g_1w,\ldots,g_dw\} \stackrel{\sim}{\ra} \psi^{-1}(\xi_0).
$$
In the sequel we identify the elements of both sets. In other words, we label the elements of the fibre
$\psi^{-1}(\xi_0)$ by $g_1w,\ldots,g_dw$.

For every point $\xi \in U$ choose a path $\gamma_{\xi}$ in $U$ connecting $\xi$ and $\xi_0$. The path
defines a bijection
$$
\mu: \psi^{-1}(\xi) \to \psi^{-1}(\xi_0)=\{g_1w=w, \dots, g_dw\}
$$
in the following way: For any $c \in \psi^{-1}(\xi)$ denote by $\tilde{\gamma}_c$ the lift of $\gamma_{\xi}$ starting at $c$.
If $g_jw \in \psi^{-1}(\xi_0)$ denotes the end point of $\tilde{\gamma}_c$, define
$$
\mu(c) = g_jw.
$$
Define
$$
K_{U,w}:=\{(x,y)\in \psi^{-1}(U) \times \psi^{-1}(U) : (\mu(x),\mu(y))-(w,w)-1 > 0\}.
$$
and let $\overline{K}_{w}$ denote the closure of $K_{U,w}$ in $C \times C$. According to
Lemma \ref{L:positivness}, $\overline{K}_{w}$ is an integral symmetric effective correspondence on $C$.\\

We claim that the correspondence $\overline{K}_{w}$ is canonically associated to the triple $(\pi, \Lambda, w)$,
i.e. does not depend on the choice of the path $\gamma_{\xi}$
connecting $\xi$ and $\xi_0$.


To see this, note that according to Remark \ref{R:sevin} the monodromy group of the covering $\psi: C \ra \P^1$
coincides with the Galois group of the Galois covering generated by $\psi$ and thus is a quotient of the Galois covering
$X$ of $\P^1$. Thus the assertion is a consequence of the $G$-invariance of the form $(\;,\;)$.

We call $\overline{K}_{w}$ the {\it Kanev correspondence} associated to the weight $w$, since it was introduced
by Kanev in \cite{kan} in a slightly different set up.
Considered as a map $C \ra \Div(C)$, it is given by


\begin{equation} \label{eq:kan}
\overline{K}_w(c)=\sum_{j=1}^d
[(\mu(c),g_jw)-(w,w)-1]\mu^{-1}(g_jw) + c.
\end{equation}

Note that $c$ is added, since in the sum $c$ appears with
coefficient $-1$, because $(\mu(c),g_iw)=(w,w)$ if $\mu(c)=g_iw$.
In the next section we need the following description of
$\overline{K}_w(c)$. Note that, according to our construction, for
any $c \in \psi^{-1}(U)$, there is a unique integer $i_c,\; 1 \leq
i_c \leq d$ such that $\mu(c) = g_{i_c}w$.

\begin{prop}\label{prop3.2}
If $c \in \psi^{-1}(U)$ with $\mu(c) = g_{i_c}w$, then
$$
\overline{K}_w(c) = \sum_{j=2}^d [(w,g_jw)-(w,w)-1]\mu^{-1}(g_{i_c}g_jw)
$$
\end{prop}

\begin{proof}


Considering Equation \ref{eq:kan} and using the $G$-invariance of
$(\;,)$,
\begin{eqnarray*}
\overline{K}_w(c) &= & \sum_{g_jw \neq \mu(c)} [(\mu(c),g_jw)-(w,w)-1]\mu^{-1}(g_jw) \\
& = & \sum_{g_jw \neq \mu(c)}[(g_{i_c}^{-1}\mu(c),g_{i_c}^{-1}g_jw)-(w,w)-1] \mu^{-1}(g_jw) \\
& =&  \sum_{g_jw \neq g_{i_c}w} [(w,g_{i_c}^{-1}g_jw)-(w,w)-1] \mu^{-1}(g_jw)\\
& = & \sum_{j=2}^d [(w,g_jw)-(w,w)-1] \mu^{-1}(g_{i_c}g_jw).
\end{eqnarray*}
\end{proof}

Next we want to work out how the Kanev correspondence $\overline{K}_w$ is related to the Schur correspondence
$\overline{S}_w$ on the curve $C$. For this we need a special choice of the set of representatives $g_1, \ldots, g_d$
of the left cosets of the subgroup $H$ of $G$.

\begin{lem} \label{lem3.3}
For any subgroup $H$ of a finite group $G$ there is a set of representatives $\{g_1, \ldots, g_d\}$ of the left cosets
$gH$ of $H$, such that $\{g_1^{-1}, \ldots, g_d^{-1}\}$ is also a set of the left cosets of $H$ in $G$.
\end{lem}

\begin{proof}
It is a consequence of the marriage theorem of combinatorics that there is a set of representatives $\{g_1, \ldots, g_d\}$
of the left cosets of $H$ in $G$ which is also a set of representatives of the right cosets of $H$ in $G$ (see e.g.
\cite{hall} Theorem 5.1.7). But if $\{g_1, \ldots, g_d\}$ are representatives of the right cosets, their inverses
$\{g_1^{-1}, \ldots, g_d^{-1}\}$ are representatives of the left cosets.
\end{proof}

In the sequel we choose a set of representatives $\{g_1, \ldots,g_d\}$ as in Lemma \ref{lem3.3}. Clearly we can
assume moreover $g_1 = 1$, implying $g_1w=w$.

\begin{prop}\label{prop3.4}
If $c \in \psi^{-1}(U)$ with $\mu(c) = g_{i_c}w$, then
$$
\overline{S}_w(c)=|H|^2\sum_{j=1}^d(w,g_jw)\mu^{-1}(g_{i_c}g_jw)$$
\end{prop}

\begin{proof}
In the proof of Proposition \ref{P:schur}, we saw that
$$
\overline{S}_w(c)=|H|\sum_{j=1}^d(w,g_jw)\sum_{h \in H}\varphi(hg_j^{-1}(x))
$$
where $x$ is a point in the fiber $\varphi^{-1}(c) \subset X$.\\
Suppose first that $\mu(c)=w$, i.e. $i_c = 1$.
As $G$ acts by left multiplication on the set $\{g_1,\ldots,g_d\}$, we have\
$$
g\mu(\varphi(x)) = \mu(\varphi(gx))
$$
for all $g \in G$.
Since $h\mu(c)=w$ for all $h\in H$, this implies
$$
\varphi(hg_j^{-1}x)=\varphi(g_j^{-1}x)=\mu^{-1}(g_j^{-1}\mu(\varphi(x))=\mu^{-1}(g_j^{-1}w)
$$
for all $j$ and $h \in H$. Therefore we get, using 
$(w,g_jw) = (g_jw,w) = (g_j^{-1}g_jw,g_j^{-1}w)$ by 
the $G$-invariance and symmetry of $(\;,\,)$,
\begin{eqnarray*}
\overline{S}_w(c) &=&|H|^2\sum_{j=1}^d(w,g_jw)\mu^{-1}(g_j^{-1} w) \\
&=&|H|^2\sum_{j=1}^d(w,g_j^{-1}w)\mu^{-1}(g_j^{-1} w)\\
&=&|H|^2\sum_{j=1}^d(w,g_jw)\mu^{-1}(g_j w),
\end{eqnarray*}
where for the last equation we used that the set of representatives $\{g_1, \ldots,g_d\}$
has the property of Lemma \ref{lem3.3}.\\
Finally, if $\mu(c) \neq w$, we have $g_{i_c}^{-1}\mu(c)=w$. Hence we can apply the above equation to the bijection
$\widetilde{\mu} = g_{i_c}^{-1} \cdot \mu: \psi^{-1}(\psi(c)) \ra \{g_1w, \ldots,g_dw\}$, which gives
\begin{eqnarray*}
\overline{S}_w(c)&=&|H|^2\sum_{j=1}^d(w,g_jw)\widetilde{\mu}^{-1}(g_j w)\\
&=&|H|^2\sum_{j=1}^d(w,g_jw)\mu^{-1}(g_{i_c}g_j w)
\end{eqnarray*}
by the definition of $\tilde{\mu}$.
\end{proof}
With the notation and identifications of above we can state the main result of this section
\begin{theorem}\label{kanev-schur}
The Kanev and Schur correspondence $\overline{K}_w$ and $\overline{S}_w$ on the curve
$C=X/H$ associated to $w$ are related as follows
$$
\overline{S}_w + |H|^2 (\Delta-\overline{K}_w)=[(w,w)+1]|H|^2 T,
$$
where $\Delta$ denotes the diagonal in $C \times C$ and $T = \psi^*\psi$ the trace correpondence of the morphism $\psi$.
\end{theorem}

\begin{proof}
It suffices to show that $[\overline{S}_w + |H|^2 (\Delta-\overline{K}_w)](c)=[(w,w)+1]|H|^2 T(c)$ for all
$c \in \psi^{-1}(U)$. But applying Propositions \ref{prop3.2} and \ref{prop3.4} we have
\begin{eqnarray*}
[\overline{S}_w + |H|^2 (\Delta-\overline{K}_w)](c) &=& |H|^2[\sum_{j=1}^d(w,g_jw)\mu^{-1}(g_{i_c}g_jw) +c \\
&& \hspace{2cm}- \sum_{j=2}^d [(w,g_jw)-(w,w)-1]\mu^{-1}(g_{i_c}g_jw)]\\
&=& |H|^2((w,w) +1) \sum_{j=1}^d \mu^{-1}(g_{i_c}g_jw)
\end{eqnarray*}
But the left multiplication with $g_{i_c}$ permutes only the elements $g_1w, \ldots, g_dw$, which implies that
$\sum_{j=1}^d \mu^{-1}(g_{i_c}g_jw) = \psi^{-1}\psi(c) = T(c)$. This completes the proof of the theorem.
\end{proof}

To every correspondence with rational coefficients on $C$ one can associate an element of the endomorphism algebra
$\End_{\Q}(JC)$ of the Jacobian $JC$ in a natural way.
If $\overline{\sigma}_w$ and $\overline{\kappa}_w$ denote the elements of $\End_{\Q}(JC)$ associated to
the correspondences $\overline{S}_w$ and $\overline{K}_w$, we get as a consequence

\begin{cor} \label{cor3.6}
The elements $\overline{\sigma}_w$ and $\overline{\kappa}_w$ of $\End_{\Q}(JC)$ associated to the Schur and
Kanev correspondences satisfy the following equation
$$
\overline{\kappa}_w = 1_{JC} + \frac{\overline{\sigma}_w}{|H|^2}.
$$
\end{cor}

\begin{proof}
This follows from the fact that the natural map from the ring of
correspondences to $\End_{\Q}(JC)$ is a homomorphism, under which
the trace correspondence maps to zero and $\Delta_C$ to $1_{JC}$.
\end{proof}

Another consequence of Theorem \ref{kanev-schur} is that $\overline{S}_w$ is a correspondence with
integer coefficients. This implies that $|H|^2 \cdot \overline{\sigma}_w$ is an endomorphism of $JC$.
The abelian subvariety
$$
B_w := Im(|H|^2 \cdot \overline{\sigma}_w)
$$
is canonically associated to the pair $(\Lambda,w)$. It certainly depends on $w$. However according to \cite{lr1} 
all abelian subvarieties $B_w$ are isogenous to each other, as long as $w \neq 0$. We call any abelian 
subvariety in the isogeny class of $B_w$ the isogeny component {\it associated to the representation} $W$.

On the other
hand, if the correspondence $\overline{K}_w$ is fixed point free,
Kanev's theorem (see \cite{lb}, Theorem 12.9.1) implies that the abelian subvariety
$$
P(C,\overline{\kappa}_w) := Im(\overline{\kappa}_w -1_{JC})
$$
is a Prym-Tyurin variety, i.e. an abelian subvariety of $JC$ such that the canonical principal polarization of $JC$
restricts to a multiple of a principal polarization. In any case, the following corollary is an immediate consequence of 
Corollary \ref{cor3.6}

\begin{cor} \label{cor3.7}  \hspace*{3.5cm} $P(C,\overline{\kappa}_w) = B_w$.
\end{cor}

Finally Theorem \ref{kanev-schur} allows us to compute the degree and the exponent of the Kanev correspondence 
$\overline{K}_w$. Recall that $d = \deg \psi = [G:H]$.

\begin{cor} \label{cor3.8} \hspace*{2.7cm} $\deg \overline{K}_w = 1 - d((w,w) +1)$. 
\end{cor}

\begin{proof} 
As we noted in Section 2, $\deg S_w = 0$ which implies $\deg \overline{S}_w = 0$. On the other hand, $\deg T = d$ and 
$\deg \Delta = 1$. So Theorem \ref{kanev-schur} gives the assertion. 
\end{proof}
Recall that the {\it exponent} of the correspondence $\overline{K}_w$ is by definition the exponent of the abelian subvariety
of $JC$ associated to it. It is the number $e$ such that $\overline{\kappa}_w$  satisfies the equation 
$\overline{\kappa}_w^2 + (e-2)\overline{\kappa}_w - (e-1) = 0$.

\begin{cor} \label{cor3.9}
The exponent of the correspondence $\overline{K}_w$ is
$$
e = e(\overline{K}_w) = - \frac{d \cdot (w,w)}{n}.
$$
\end{cor}
Note that $(w,w)$ is negative. So $e$ is a positive integer.
\begin{proof}
Since $p_w$ is an idempotent, we get from the definition of $S_w$ in Section 2, 
$S_w^2 = \frac{|G| \cdot (w,w)}{n}S_w.$
This means for the associated endomorphism $\sigma_w$ of $JX$
$$
\sigma_w^2 =  \frac{|G| \cdot (w,w)}{n}\sigma_w.
$$
If $\varphi_* : JX \ra JC$ denotes the norm map and $\varphi^* : JC \ra JX$ 
the pull-back map associated to the morphism $\varphi: X \ra C$, the endomorphisms $\sigma_w \in \End(JX)$ and 
$\overline{\sigma}_w \in \End(JX)$ are related by
$$
\overline{\sigma}_w = \varphi_* \sigma_w \varphi^*.
$$
This implies, since $\varphi^* \varphi_*$ is multiplication by $\deg \varphi = |H|$ on Im $\varphi^* \subset JX$, 
$$
\overline{\sigma}_w^2 = \varphi_* \sigma_w \varphi^* \varphi_* \sigma_w \varphi^*
 = |H| \cdot \varphi_* \sigma_w^2 \varphi^*
=  \frac{ |H| |G| (w,w)}{n} \cdot \overline{\sigma}_w.
$$
Applying Corollary \ref{cor3.6}, this gives
\begin{equation} \label{eq3.4}
(1_{JC} - \overline{\kappa}_w)^2 = \frac{\overline{\sigma}_w^2}{|H|^4} 
= \frac{|G| (w,w)}{|H| \cdot n} \cdot \frac{\overline{\sigma}_w}{|H|^2} = e \cdot (1_{JC} - \overline{\kappa}_w)
\end{equation}
with $e = - \frac{d \cdot (w,w)}{n}$ or equivalently 
$\overline{\kappa}_w^2 +(e-2) \overline{\kappa}_w - (e-1) = 0$.
\end{proof}

\begin{rem}
In the set up of \cite{kan} Corollaries \ref{cor3.8} and \ref{cor3.9} have been proven in Propositions 5.2 and 5.3.
\end{rem}

\section{Dimension of $P(C, \kappa_w)$ and fixed points of $\overline{K}_w$}\label{S:fixed-kanev}

\noindent
Let the notation be as above. In particular $X$ is a curve with
$G$-action, $w$ a weight of the representation $\Lambda$, and we have diagram  (\ref{eq:diagram}). First we compute
the dimension of the abelian subvariety $P(C,\overline{\kappa}_w)$ associated to the Kanev correspondence
$\overline{K}_w$. Note that $e$ and $\deg \overline{K}_w$ have been computed in Corollaries \ref{cor3.8} and \ref{cor3.9}.

\begin{prop} \label{prop4.1} \hspace{1.5cm}
$\dim P(C,\overline{\kappa}_w) = \displaystyle\frac{1}{e}(g_C - \deg \overline{K}_w + \frac{(\overline{K}_w \cdot \Delta)}{2}).$
\end{prop}

\begin{proof}
According to (\ref{eq3.4}) the element $\frac{1}{e}(1_{JC} - \overline{\kappa}_w) \in \End_{\Q}(JC)$ is the symmetric idempotent corresponding to
the abelian subvariety $P(C,\overline{\kappa}_w)$. Hence \cite{lb} Corollary 5.3.10 gives
$$
\dim P(C, \overline{\kappa}_w) = Tr_a( \frac{1}{e}(1_{JC} - \overline{\kappa}_w)) 
= \frac{1}{e}(g_C - \frac{1}{2} Tr_r(\overline{\kappa}_w)).
$$
On the other hand, according to \cite{lb} Proposition 11.5.2 we have for the rational trace of $\overline{\kappa}_w$
$$
Tr_r(\overline{\kappa}_w) = 2 \deg \overline{K}_w - (\overline{K}_w \cdot \Delta)
$$
which completes the proof of the proposition.
\end{proof}

In order to compute the number $(\overline{K}_w \cdot \Delta)$ of fixed points of the Kanev correspondence,
we need to know more about the ramification of the map $\pi:X \ra \P^1$. 
Let $B = \{b_1,\dots, b_s\} \in \P^1$ denote the branch locus of $\pi:X \ra \P^1$ which we assume also to be the branch locus 
of $\psi: C \ra \P^1$. According to Remark \ref{R:sevin} this assumption is fullfilled if the represesentation $W$
is faithful, which is the case for all interesting examples. Denote $U=\P^1 \setminus B$ and 
$$
m:\pi_1(U,\xi_0) \to G \in S_d
$$ the
monodromy homomorphism of $\psi$ with respect to some base point $\xi_0 \in U$.

Let $\{C_1, \dots, C_r\}$ denote the set of non-trivial conjugacy classes
of cyclic subgroups of $G$. We say that the branch point $b_i \in B$ of the map $\pi$ is \textit{of type $C_j$} 
if the following condition is satisfied:
Let $\gamma_i$ be a standard loop composed by a path $\gamma$ starting at $\xi_0$, a
small positively oriented circle around $b_i$ and $\gamma^{-1}$. Then $m(\gamma_i)$ is contained in some representative of the class $C_i$.
In other words, considering the action of $G$ on $X$, the stabilizers of the points
in the fibre of $b_i \in B$ belong to the conjugacy class $C_j$.\\

We assume in the sequel that $\alpha_j$ branch points  $b_i \in B$ are of type $C_j$ for $j=1, \ldots r$. In particular
$\sum_{j=1}^r \alpha_j = s$. In terms of these data the genus of $C$ has been computed in \cite{ksir}, equation $(10)$ 
and independently in \cite{yo}, Corollary 3.4. If $G_j$ denotes a representative of the class $C_j$, 
\begin{equation} \label{eq4.1}
g_C = -d +1 + \frac{1}{2} \sum_{j=1}^r \alpha_j(d - |H \backslash G /G_j|)
\end{equation}
where $|H \backslash G /G_j|$ denotes the number of double cosets of the subgroups $H$ and $G_j$ in $G$.
Moreover, using Corollary \ref{cor3.7}, \cite{yo}, Corollary 5.12 gives
\begin{equation} \label{eq4.2}
\dim P(C,\overline{\kappa}_w) = -n + \frac{1}{2} \sum_{j=1}^r \alpha_j ( n - \dim Fix_{G_j}(W)).
\end{equation}
Applying these results we obtain with the notation as above for the number of fixed points $(\overline{K}_w \cdot \Delta)$
of $\overline{K}_w$,

\begin{prop} \label{prop4.2} 
$(\overline{K}_w \cdot \Delta) = \sum_{j=1}^r \alpha_j \left(e(\dim W-\dim Fix_{G_j}(W)) - ([G:H]- |H \backslash G /G_j|)\right).$
\end{prop}

\begin{proof}
Applying Proposition \ref{prop4.1} we have, using (\ref{eq4.1}) and (\ref{eq4.2})
\begin{eqnarray*}
(\overline{K}_w \cdot \Delta) &=& 2(e \dim P(C, \overline{\kappa}_w) +\deg \overline{K}_w -g_C) \\
&=& \sum_{j=1}^r \alpha_j (en - d + |H \backslash G /G_j| - e \dim Fix_{G_j} (W)) + 2(\deg \overline{K}_w +d - 1 - en).  
\end{eqnarray*}
But $\deg \overline{K}_w +d - 1 - en = 0$ according to Corollaries \ref{cor3.8} and \ref{cor3.9}. This implies the assertion 
recalling that $n = \dim W$ and $d = [G:H]$.
\end{proof}

\begin{rem}
Recall that $[G:H]=\dim \mbox{Ind}^G_H$ and $|H \backslash G /G_j|=\dim Fix_{G_j}(\mbox{Ind}^G_H)$, where $\mbox{Ind}^G_H$ is the 
representation of $G$ induced by the trivial representation of $H$. Thus the fixed points of the correspondence are contained 
in the ramification points of the covering and such a point is a fixed point of $\overline{K}_w$ if and only if the codimension 
of $Fix_{G_j}(\mbox{Ind}^G_H)$ in $\mbox{Ind}^G_H$ is strictly smaller than the $e$-fold of the codimension of $Fix_{G_j}(W)$ in $W$.

\end{rem}



\section{Examples}

\subsection{The Weyl group of $E_6$}

Let $G = W(E_6)$ denote the Weyl group of type $E_6$. We use the notation of Bourbaki \cite{bou}. 
Let $\Lambda$ denote the root lattice. So $W = \Lambda \otimes \Q$ is the root representation and the action 
of $G$ is generated by the 6 reflections $s_{\alpha_i}$ associated to the simple roots $\alpha_i$ for $i = 1, \ldots, 6$. 
As a weight we choose the fundamental weight 
$\varpi_6 = \frac{1}{3}(2\alpha_1 + 3 \alpha_2 + 4\alpha_3 + 6\alpha_4 + 5\alpha_5 + 4\alpha_6)$. With respect to the basis 
$\alpha_1, \ldots, \alpha_6$  bilinear form $(\;,\;)$ is given by the negative of the Cartan matrix. This implies 
$(\varpi_6,\varpi_6) = - \frac{4}{3}$. The stabilizer $H$ of $\varpi_6$
is the subgroup of $G$ generated by $\{s_{\alpha_i}\;|\; \; 1 \leq i \leq 5\}$. It is of index $d = 27$ in $G$. 
Finally, let $X$ denote a curve with $G$-action, with quotient $X/G = \P^1$ and associated diagram \eqref{eq:diagram}.\\ 

Choosing a set of representatives $g_i, \; 1 \leq i \leq 27$ for the left cosets of $H$ in $G$ we have 
with the notation of Proposition \ref{prop3.2} for the associated Kanev correspondence
$$
{\overline K}_{\varpi_6}(c) = \sum_{j=1}^{27} a_{ij}\mu^{-1}(g_j\varpi_6),
$$
where $i$ runs from 1 to 27, representing the case $\mu(c) = g_i \varpi_6$. Computing the matrix $(a_{ij})$ (using Maple) 
one obtains the incidence matrix of the 27 lines on a general cubic surface as it should be (see \cite{kan}).

The group $G=W(E_6)$ has $24$ non-trivial cyclic subgroups up to conjugacy. We distinguish two of them:
\begin{itemize}
\item $C_{1}$, the class of subgroups of relections i.e. conjugated to $G_{1}=<s_{r_1}>$,\\ 
\item $C_{2}$, the class of cyclic subgroups of order $3$ conjugated to $G_{2}=<s_{r_1}s_{r_2}>$.
\end{itemize}
Note that the points of type $C_{1}$ correspond to points of ''simplest
ramification'', whereas
the points of type $C_{2}$ are not allowed in \cite{kan}.

Using Corollaries \ref{cor3.8} and \ref{cor3.9} and Proposition \ref{prop4.2} and noting that 

\vskip12pt

\begin{tabular}{l l l l}

$|H \backslash G/G_{1}| = 21$, & $|H \backslash G/G_{2}| = 15$, & $\dim Fix_{G_1}W=5$, & $\dim Fix_{G_2} W=4$
\end{tabular}

\vskip12pt

\noindent we conclude

\begin{prop} \label{lem5.1}
Suppose the Galois covering $\pi:X \ra \P^1$ admits only branch points of type $C_{1}$ and $C_{2}$. 
Then the Kanev correspondence
${\overline K}_{\varpi_6}$ is fixed point-free with exponent 6 and degree 10. 
\end{prop}

Notice that with Proposition \ref{prop4.2} one can also show that any other branch point leads to fixed 
points of ${\overline K}_{\varpi_6}$.

\begin{cor} \label{cor5.2}
Assume in addition that $\pi:X \ra \P^1$ admits $\alpha$ branch points of type $C_{1}$, $\beta$ 
branch point of type $C_{2}$ and no branch points of another type.
Then the abelian variety $P(C,\kappa_{\varpi_6})$ is a Prym-Tyurin variety of exponent 6 and
dimension $\frac{\alpha}{2} + \beta - 6$ in the Jacobian of $C$ of dimension $g(C) = 3\alpha + 6 \beta - 26$.   
\end{cor} 

\begin{proof}
According to Proposition \ref{lem5.1} the correspondence ${\overline K}_{\varpi_6}$ fixed point-free of exponent 6. Hence 
Kanev's criterion (see \cite[Theorem 12.9.1]{kan}) implies the first assertion. The formula for the dimension follows from 
Proposition \ref{prop4.1} and equations \eqref{eq4.1} and \eqref{eq4.2}.
\end{proof}

We are left with the question of existence of such a Galois covering. 

\begin{prop} \label{prop5.3}
For $\alpha = 12 + 2k$ and $\beta = 2n$ with non-negative integers $k$ and $n$ there exist Galois coverings $\pi: X \ra \P^1$
with group $G = W(E_6)$ ramified exactly over $\alpha$ points of type $C_1$ and $\beta$ points of type $C_2$. The dimension 
of the corresponding Prym-Tyurin variety is in this case 
$$
\dim P(C, \kappa_{\varpi_6}) = k + 2n.
$$ 
\end{prop}

\begin{proof}
According to \cite{yo}, the question of existence is equivalent to finding
a subset of $C_1$ or $C_2$ or both, generating $G$ and having trivial product.
Since $G$ is generated by $6$ reflections, we can certainly generate $G$ with $12$ reflections with trivial product.
On the other hand, it is not possible to generate $G$ with elements $C_{2}$ alone. However we can add to the $12$ reflections
an even number of elements of $C_1$ and $C_2$ (add with an element also its inverse) keeping the trivial product.
This proves the first assertion. The last assertion follows from Corollary \ref{cor5.2}.
\end{proof}

\begin{rem}
\label{rem5.4}
The coverings with no branch point of type $C_2$ are exactly the coverings ''with simplest ramification'' occurring in 
\cite{kan}. Hence our construction gives new examples of Prym-Tyurin varieties.
\end{rem}

\begin{rem}
The family $\mathcal{R}_G(\alpha,\beta)$ of Galois coverings as in Proposition \ref{prop5.3} 
is of dimension $\alpha+\beta-3=9+2n+2k$. On the other hand, the moduli space $\A_{k+2n}$ 
of principally polarized abelian varieties
of dimension $k+2n$ has dimension
$\frac{(k+2n)(k+2n+1)}{2}$.

Hence $\dim \mathcal{R}_G(24,0) = \dim \A_6 = 21$ and similarly  $\dim \mathcal{R}_G(14,2) = \dim \A_5 = 15$. 
The natural question to ask is of course: Is a general principally polarized abelian variety of dimension 6 
a Prym-Tyurin variety of exponent 6 in this way?                                                                                     
\end{rem}

\subsection{The alternating group} Let $G= A_n$ denote the alternating group of degree $n$ for $n \geq 4$. Denote by 
$e_1, \ldots, e_n$ the standard basis of $\R^n$. Then 
$$
\Lambda = \; < E_1= e_1-e_2, \ldots, E_{n-1} = e_{n-1}-e_n >
$$
is an absolutely irreducible $\Z$-representation of rank $n-1$ of $G$. As a weight we choose 
$$
w=\sum_{i=1}^{n-2}e_i-\frac{n-2}{n}\sum_{i=1}^n e_i=\frac{1}{n}\left[\sum_{i=1}^{n-2}2iE_i+(n-2)E_{n-1}\right]
$$ 

The bilinear form $(\;,\;)$ is given by
$$
(E_i,E_j) = \left\{\begin{array}{ccc}
-2 &  & i=j,\\
1 & for & i=j+1 \; \mbox{or} \; i=j-1,\\
0 && \mbox{otherwise}.
\end{array}\right.
$$
This implies 
$$
(\omega,\omega) = -\frac{2(n-2)}{n}.
$$
Clearly the stabilizer of $\omega$ in $A_n$ is $H = M \cap A_n$ where $M =  S_{n-2}^{\{n-1,n\}}\times <(n-1\;\;\; n)>$ 
and we have
$$
[A_n:H] = \frac{n(n-1)}{2}.
$$
Notice that with the notation of \cite{bou}, $\Lambda$ is the restriction of the root lattice of the symmetric group 
of degree $n$, the weight is $\omega = \varpi_{n-2}$ and the bilinear form is given by the negative of the Cartan matrix.

Finally let $X$ denote a curve with $G$-action, with quotient $X/G = \P^1$ and associated diagram \eqref{eq:diagram}.
Let $C_3$ denote the conjugacy class of cyclic subgroups of order 3 generated by elements of the form $(ijk)$. 
\begin{prop}
Suppose the Galois covering $\pi:X \ra \P^1$ admits only branch points of type $C_3$.
Then the Kanev correspondence ${\overline K}_\omega$ is fixed point-free with exponent $n-2$ and degree $\frac{1}{2}(n^2-5n+6)$.
\end{prop}

\begin{proof}
Let $G_3 \in C_3$ any representative of the class. According to Proposition \ref{prop4.2} we have to compute $\dim \Fix_{G_3}W$
and $|H\bs G/G_3|$.
First note that
$$
\dim \Fix_{G_3} W=\;<\Ind_{G_3}^{A_n}1_{G_3}, \Lambda>,
$$
\noindent 
where $<\;,\;>$ denotes the character product of 
representations and $\Ind_{G_3}^{A_n}1_{G_3}$ is the induced representation from the trivial representation of $G_3$ to $A_n$.
By Frobenius reciprocity we conclude
$$
\dim Fix_{G_3} W= \; <1_{G_3}, Res_{G_3} \Lambda>\; 
=\frac{1}{|G_3|}\sum_{g\in G_3}\chi_{1_{G_3}}(g^{-1})\chi_{\Lambda}(g)=n-3.
$$
For the double cosets cardinality $|H\bs G/G_3|$, note that $H\bs A_n$ is in bijective correspondence 
with the orbit of $\{n-1,n\}$ in the set $P$ 
of unordered pairs and that $\sigma \in A_n$ acts from the left on $P$ mapping $\{a,b\}$ to $\{\sigma(a),\sigma(b)\}$.
Hence we just need to count the orbits $P/G_3$ of $G_3$ in the set $P$. 


Now we use Polya-Burnside formula:
$$
|P/G_3|=(1/|G_3|)\sum_{g\in G_3} |P^g|.
$$
\noindent Here $G_3=\{id,g_1,g_1^2\}$, thus for $g=id$ we have $|P^g|=\binom{n}{2}$, because it fixes all the pairs.
For the other two elements of $G_3$ the cardinality of $P^g$ is $\binom{n-3}{2}$, 
because an element $(i\; j\; k) \in G_3$ fixes exactly the pairs not intersecting the set $\{i,j,k\}$.

We conclude that the number of orbits is $\frac{1}{3}(\binom{n}{2}+ 2\binom{n-3}{2})$, which gives

$$
|H\bs G/G_3|=\frac{n^2-5n+8}{2}.
$$

Using this we obtain from Proposition \ref{prop4.2} that $\overline{K}_w$ is fixed point-free.
Finally Corollaries \ref{cor3.8} and \ref{cor3.9} give the formulas for $e$ and $\deg \overline{K}_w$.
\end{proof}

The proof of the following corollary is the same as for Corollary \ref{cor5.2}.

\begin{cor}
Assume in addition that $\pi: X \ra \P^1$ admits $\alpha$ branch points of type $C_3$ and is brancj point-free elsewhere. 
Then the abelian variety $P(C,\omega)$ is a Prym-Tyurin variety of exponent $n-2$ and dimension $\alpha -n + 1$ in the  
Jacobian of $C$ of dimension $\alpha(n-2) + 1 - \binom{n}{2}$.
\end{cor}

Again we have to check the existence of such Galois coverings. In the same way as Proposition \ref{prop5.3} we deduce

\begin{prop}
For $\alpha \geq 2(n-2), \; n \geq 4$ there exist Galois coverings with group $A_n$ ramified exactly over $\alpha$
points of type $C_3$.
\end{prop}

\begin{rem}
The abelian variety $P(C, \omega)$ is isogenous to a Jacobian variety. In fact, consider the subgroup $N=A_n\cap S_{n-1}$
of index $n$ in $A_n$, and the quotient curve $X/N$. It is easy to see using \cite{yo} its genus is given by
$$
g_{X/N}=-|A_n:N|+1+\frac{1}{2}\alpha (|A_n:N|-|N\bs A_n / G_3|)=\alpha-n+1\;.
$$
Moreover, the representation of $A_n$ induced by the trivial one of $N$ decomposes as $1 \oplus \Lambda$
therefore the Jacobian variety $JX/N$ of the curve $X/N$ is isogenous to the factor $B_\omega$ corresponding to
$\Lambda$ in the isotypical decomposition of $JX$ (see \cite{lr1}). But we showed that the Schur and the Kanev correspondence 
yield the same abelian subvariety, implying that $JX_N$ is isogenous to $P(C,\omega)$.
\end{rem}

\end{document}